\documentclass[a4paper, 11pt]{article}
\usepackage{amsmath, amsfonts,amssymb}
\usepackage[english, russian]{babel}
\begin {document}

\begin{center}
\textbf {Geometrical properties of the space $I_f(X)$ \\ of idempotent probability measures} \\

\medskip
\textbf {Zaitov A.\,A.$^{1}$, Ishmetov A.\,Ya.$^{2}$}\\
\smallskip
{$^{1}$ Tashkent institute of architecture and civil engineering}\\
{adilbek\_zaitov@mail.ru},\\
{$^{2}$ Tashkent institute of architecture and civil engineering}\\
{ishmetov\_azadbek@mail.ru}\\
\end{center}

\begin{abstract}
In this paper we prove that for a compact space $X$ inclusion $I_{f}(X)\in ANR$ holds if and only if $X\in ANR$. Further, it is shown that the functor $I_{f}$ preserves property of a compact to be  $Q$-manifold or a Hilbert cube, properties of maps fibres to be $ANR$-compact, $Q$-manifold, Hilbert cube (the finite of Hilbert cube).\\

\end{abstract}

The theory of idempotent measures belongs to idempotent mathematics, i. e. the field of the mathematics based on replacement of usual arithmetic operations with idempotent (as, for example, $x\oplus y=\max\{x,y\}$).  The idempotent mathematics intensively develops at this time (see, for example, [1], survey article [2] and the bibliography in it). Its communication with traditional mathematics is described by the informal principle [2] according to which there is a heuristic compliance between important, interesting and useful designs
the last and similar results of idempotent mathematics.

In the present article we investigate a subfunctor $I_f$ of a functor of idempotent probability measures in category compact Hausdorff spaces. In traditional mathematics to it there corresponds the functor $P_f$ of probability measures. The concept of an idempotent measure (Maslov's measure) finds numerous applications in various field of mathematics, mathematical physics and economy. In particular, such measures arise in problems of dynamic optimization [3]; the analogy between Maslov's integration and optimization is noted also in [4]. In [5] it is claimed that use of measures of Maslov for modeling of uncertainty in mathematical economy can be so relevant as far as also use of classical probability theory.

Unlike a case of probability measures to which consideration extensive literature is devoted (see [6]) geometrical and topological properties of the spaces of idempotent measures were practically not investigated.\\

Consider set $\mathbb{R}_{+}=\mathbb{R}\bigcup \left\{ -\infty  \right\}$  with two algebraic operations: addition $\oplus $ and
multiplication $\odot$ determined as follows: $u\oplus v=\max \{u,v\}$ and $u\odot v=u+v$, $u,v\in \mathbb{R}_{+}$, where $\mathbb{R}$ is the set of real numbers.

Let $X$ be a compact Hausdorff space ($\equiv$ a compact), $C(X)$ be the algebra of continuous functions on $X$ with usual algebraic
operations. On $C(X)$ operations $\oplus $  and $\odot $ we will determine by rules $\varphi \oplus \psi =\max \{\varphi ,\psi \}$ and $\varphi \odot \psi=\varphi + \psi$ where $\varphi $, $\psi \in C(X)$.

Remind a functional $\mu :C(X)\to \mathbb{R}$ is called $\left[ 7 \right]$ to be an idempotent probability measure on $X$ if
it satisfies the following properties:

(1) $\mu ({{\lambda }_{X}})=\lambda $ for all $\lambda \in \mathbb{R}$, where  ${{\lambda }_{X}}$ -- constant function;

(2) $\mu(\lambda \odot \varphi)=\lambda \odot \mu(\varphi)$ for all $\lambda \in \mathbb{R}$ и $\varphi \in C(X)$;

(3) $\mu (\varphi \oplus \psi )=\mu (\varphi )\oplus \mu (\psi )$ for all $\varphi $, $\psi \in C(X)$.

For a compact $X$ we denote by $I(X)$ the set of all idempotent probability measures on $X$.
We consider $I(X)$ as subspace of ${{\mathbb{R}}^{C(X)}}$.

For given compacts $X$, $Y$ and  continuous map $f:X\to Y$  it is possible to check that determined by a formula $I(f)(\mu )(\psi )=\mu (\psi \circ f)$ naturally way arised map $I(f):I(X)\to I(Y)$ continuous. The construction  $I$ is a normal functor. Therefore  for each idempotent probability measure  $\mu \in I(X)$ one may determine its support:

$$\text{supp}\mu =\bigcap \left\{ A\subset X:\overline{A}=A,\,\,\mu \in I(A) \right\}.$$

For a positive integer $n$ and compact $X$ define the following set
$${{I}_{n}}(X)=\left\{ \mu \in I(X):\,\,\left| \operatorname{supp}\mu  \right|\le n \right\}.$$

Put
$${{I}_{\omega }}(X)=\bigcup\limits_{n=1}^{\infty }{{{I}_{n}}(X)}.$$

Set ${{I}_{\omega }}(X)$ everywhere dense $\left[ 7 \right]$ in $I(X)$. Idempotent probability measure
$\mu \in {{I}_{\omega}}(X)$  is called an idempotent probability measure with finite support.

We determine a subfunctor $I_{f}$ of the functor $I$ of idempotent probability measures. For a given compact $X$ a set $I_{f}(X)$ consists of measures with finite support such that if support of a measure $\mu$ consists from  $n$  points $x_{1}, x_{2},..., x_{n}$, that only a unique of these points $max$-$plus$-barycenter mass is equal to zero, all other masses are not big $-ln(n+1)$. By definition
 $$I_f(X)=\{\mu=\bigoplus_{i=1}^n\lambda_i\odot\delta_{x_i} \in I_\omega(X):\mbox{ equality } \lambda_{i_0}=0 \mbox{ holds for only unique }$$
 $$\mbox { index } i_0 \mbox{ and } \lambda_i\leq -ln(n+1) \mbox{ for all } i\in \{1, ..., n\}\setminus \{i_0\} \}$$

This functor is interesting that it is a functor with the finite support, and has not finite degree. Functor $I_{f}:Comp\rightarrow Comp$  is a normal subfunctor of a functor $I$ idempotent probability measures.

From definition of elements of the space $I_{f}(X)$ follows that a set $\delta(X)$  Dirac's measures lies in $I_{f}(X)$.

Let $X$ и $Y$ be two compacts lying in spaces $M$ and $N$  respectively, where  $M,N\in AR$. A sequence of maps $f_{k}:M\rightarrow N$, $k=1,2,...,$ is called to be fundamental sequence from $X$ into $Y$, if for each neighbourhood  $V$ of the compact $Y$  (in $N$) there is such neighbourhood $U$ of the compact $X$   (in $M$), that $f_{k}|_{U}=f_{k+1}|_{U}$  at $V$ almost for all $k$. It means that there is such homotopy $f_{k}:U\times[0,1]\rightarrow V$, что $f_{k}(x,0)=f_{k}(x)$ и $f_{k+1}(x,1)=f_{k+1}(x)$  for all $x\in U$. We will denote this fundamental sequence through $\left\{f_{k},X,Y\right\}$  or shortly through $f$, also we will write $f:X\rightarrow Y$.

A fundamental sequence $f=\left\{f_{k},X,Y\right\}$ is generated by map $f:X\rightarrow Y$ if  $f_{k}(x)=f(x)$  for all $x\in X$ and for all $k=1,2,...$.

Spaces $X$ and $Y$ are fundamentally equivalent if there are such two fundamental sequences $f:X\rightarrow Y$ and $g:Y\rightarrow X$ that $gf=id_{X}$ and $fg=id_{Y}$.

The relation of fundamental equivalence is the equivalence relation therefore the class of all spaces decomposes to in pairwise disjoint classes of spaces which are called shape. So two spaces belong to the same shape if and only if when they are fundamentally equivalent. A shape containing a space $X$ is called [10] a shape of the space $X$ and denote by $Sh(X)$. It is known that for two neighbourhood retract of $A$ and $B$ fairly $Sh(A)=Sh(B)$  if and only if when they are homotopically equivalent.

Let's remind that maps  $f:X\rightarrow Y$ is called cellular and similar (briefly -- $CE$) [11], if any compact $A\subset Y$  prototype  $f^{-1}(A)$  is a compact and for each point $y\in Y$  prototype $f^{-1}(y)$   has sheyp points (i.e. a prototype $f^{-1}(y)$  homotopically it is equivalent to a point).

If $r:X\rightarrow F$ is a retraction and there also exists such a homotopy $h:X\times[0,1]\rightarrow F$ that $h(x,0)=x$, $h(x,1)=r(x)$ for all $x\in X$ then $r$ is deformation retraction, and $F$ is deformation retract of the space $X$. Deformation retraction  $r:X\rightarrow F$  is strongly deformation retraction if for a homotopy $h:X\times[0,1]\rightarrow F$ we have $h(x,t)=x$ for all $x\in F$  and all $t\in[0,1]$  [12].

Take any measure $\mu=\bigoplus\limits_{i=1}^{n}\lambda_{i}\odot\delta_{x_{i}}\in I_{f}(X)$. Let $\lambda_{i_{0}}=0>-ln(n+1)$. Measure $\mu$  puts in correspondence a point $\delta_{x_{i_{0}}}$ of a compact $\delta(X)$. The obtained correspondence $I_{f}(X)\rightarrow\delta(X)$  denote by $r_{\delta(X)}^{I(X)}$.

{\bf Теорема 1.} {\it For any compact $X$ a map $r_{\delta(X)}^{I(X)}:I_{f}(X)\rightarrow \delta(X)$ is continuous, open, cellular-similar (all fibres are collapsible) retraction.}

{\bf Proof.} By construction the map $r_{\delta(X)}^{I(X)}:I_{f}(X)\rightarrow \delta(X)$ is defined correctly. It is clear that $r_{\delta(X)}^{I(X)}(\delta_{x})=\delta_{x}$ for each  $x\in X$, i. e. every point of the space $\delta(X)$ is fixed-point according to the map  $r_{\delta(X)}^{I(X)}$. So, it is established that $r_{\delta(X)}^{I(X)}$ is retraction.

It is clear, that for each point $x\in X$ the fibre $\left(r_{\delta(X)}^{I(X)}\right)^{-1}(\delta_{x})$  is compact. It is clear, that for each point $\mu\in \left(r_{\delta(X)}^{I(X)}\right)^{-1}(\delta_{x})$ an interval $[\mu,\delta_{x}]=\{\lambda_{1}\odot\delta_{x}\oplus\lambda_{2}\odot\mu: \lambda_{1}\oplus\lambda_{2}=\mathbf{1};\,\, \lambda_{1}\geq\mathbf{0},\lambda_{2}\geq\mathbf{0}\}$  lies in the fibre $\left(r_{\delta(X)}^{I(X)}\right)^{-1}(\delta_{x})$.

Fix a fibre $\left(r_{\delta(X)}^{I(X)}\right)^{-1}(\delta_x)$ and define a map  $h:\left(r_{\delta(X)}^{I(X)}\right)^{-1}(\delta_x)\times[0,1]\rightarrow \left(r_{\delta(X)}^{I(X)}\right)^{-1}(\delta_x)$  by the rule
$$h(\mu,t)=\left(ln(1-t)-lnt\oplus ln(1-t)\right)\odot\delta_{x}\oplus\left(lnt-lnt\oplus ln(1-t)\right)\odot\mu,$$
where $\mu=\bigoplus\limits_{i=1}^{n}\lambda_i\odot\delta_{x_i}\in \left(r_{\delta(X)}^{I(X)}\right)^{-1}(\delta_{x})$ and $t\in [0,\ 1]$.

It is easy to see that  $h$  is a homotopy, connecting identity mapping,  $h_1=id_{\left(r_{\delta(X)}^{I(X)}\right)^{-1}(\delta_x)}$ and stationary map $h_0=\left(r_{\delta(X)}^{I(X)}\right)^{-1}(\delta_x)\rightarrow\{\delta_x\}$.

So each pre-image has shape of a point, i. e. a retraction $r_{\delta(X)}^{I(X)}$  is cellularity-similar. Moreover, these pre-images subtend to a point.

By virtue of $I_f(X)=\bigcup\limits_{x\in X}\left(r_{\delta(X)}^{I(X)}\right)^{-1}(\delta_x)$ and $\left(r_{\delta(X)}^{I(X)}\right)^{-1}(y)\cap\left(r_{\delta(X)}^{I(X)}\right)^{-1}(z)=\varnothing$ on $y\neq z$, we obtain that each fibre $\left(r_{\delta(X)}^{I(X)}\right)^{-1}(x)$, $x\in X$, is open and closed in $I_f(X)$. Therefore map $r_{\delta(X)}^{I(X)}:I_{f}(X)\rightarrow\delta(X)$ is continuous.

It is remain to show the map $r_{\delta(X)}^{I(X)}$ is open. Let $\mu\in I_{\omega}(X)$,  $\mu=\bigoplus\limits_{i=1}^n \lambda_{i}\odot\delta_{x_i}$.
Define a set
$$\left\langle\mu;U_{1},...,U_{n};\varepsilon\right\rangle=\left\{\nu=\bigoplus\limits_{j=1}^k \gamma_i\odot\delta_{y_i}\in I_{\omega}(X): supp\nu\cap U_i\neq\varnothing, i=1,2,...,n,\right.$$  $$\left.supp\nu\subset\bigcup\limits_{i=1}^n U_i, \ \ \bigoplus\limits_{y_j\in U_i}|\lambda_i-\gamma_j|<\varepsilon\right\},$$
where  $U_i$  are open neighbourhoods of points $x_i$, $i=1, ..., n$, respectively. It is clear, that set $\left\langle\mu;U_{1},...,U_{n};\varepsilon\right\rangle$ is open in  $I_{\omega}(X)$.
Let's show that sets of a type $\left\langle\mu;U_{1},...,U_{n};\varepsilon\right\rangle$ form base of weak topology in $I_\omega(X)$. Let $\left\langle\mu;\varphi;\varepsilon\right\rangle$ be an element of the prebase where $\varphi\in C(X)$, $\varepsilon>0$ and $\mu=\bigoplus\limits_{i=1}^{n}\lambda_{i}\odot\delta_{x_{i}}\in I_{\omega}(X)$. As function $\varphi$ is continuous, for each point $x_i$ there is its open neighbourhood $U_i$ such that for any point $y\in U_i$ the inequality $|\varphi(x_i)-\varphi(y)|<\frac{\varepsilon}{2}$ holds. Further, for every $\nu=\bigoplus\limits_{j=1}^k\gamma_j\odot\delta_{y_j}\in\left\langle\mu;U_{1},...,U_{n}; \frac{\varepsilon}{2}\right\rangle$ we have $\bigoplus\limits_{y_i\in U_i}\left|\lambda_i-\gamma_j\right|<\frac{\varepsilon}{2}$. Let's estimate an absolute value $|\mu(\varphi)-\nu(\varphi)|= \left|\bigoplus\limits_{i=1}^n\lambda_i\odot\varphi(x_i) - \bigoplus\limits_{j=1}^k\gamma_j\odot\varphi(y_j)\right|=a$. Two cases are possible:

1) $\bigoplus\limits_{i=1}^n\lambda_i\odot\varphi(x_i)\geq \bigoplus\limits_{j=1}^k\gamma_j\odot\varphi(y_j)$. Let $\bigoplus\limits_{i=1}^n\lambda_i\odot\varphi(x_i)=\lambda_{i'}\odot\varphi(x_{i'})$. Then

$a=\bigoplus\limits_{i=1}^n\lambda_i\odot\varphi(x_i)-\bigoplus\limits_{j=1}^k\gamma_j\odot\varphi(y_j) =\lambda_{i'}\odot\varphi(x_{i'})-\bigoplus\limits_{j=1}^k\gamma_j\odot\varphi(y_j)\leq(\mbox{for each } y_{j}\in U_{i'})$. $\leq\lambda_{i'}\odot\varphi(x_{i'})-\gamma_{j}\odot\varphi(y_{j})= |\lambda_{i'}\odot\varphi(x_{i'})-\gamma_{j}\odot\varphi(y_{j})|\leq|\lambda_{i'}-\gamma_{j}|+ |\varphi(x_{i'})-\varphi(y_{j})|<\varepsilon$.

2) $\bigoplus\limits_{i=1}^n\lambda_i\odot\varphi(x_i)\leq\bigoplus\limits_{j=1}^k\gamma_j\odot\varphi(y_j)$. Let $\bigoplus\limits_{j=1}^n\gamma_j\odot\varphi(y_j)=\gamma_{j'}\odot\varphi(y_{j'})$. Then

$a=\bigoplus\limits_{j=1}^n\gamma_j\odot\varphi(y_j)-\bigoplus\limits_{i=1}^n\lambda_i\odot\varphi(x_i)= \gamma_{j'}\odot\varphi(y_{i'})-\bigoplus\limits_{i=1}^n\lambda_i\odot\varphi(x_i)\leq(\mbox{for all } {i}, \mbox{ for which } y_{j'}\in U_{i}) \leq\gamma_{i}\odot\varphi(y_{j'})-\lambda_{i}\odot\varphi(x_{i})= |\gamma_{j'}\odot\varphi(y_{j'})-\lambda_{i}\odot\varphi(x_{i})|\leq|\lambda_{i}- \gamma_{j'}|+|\varphi(x_{i})-\varphi(y_{j'})|<\varepsilon$.

So, $a<\varepsilon$ for two cases, i. e. $|\mu(\varphi)-\nu(\varphi)|<\varepsilon$. From here $\nu\in\left\langle\mu;\varphi;\varepsilon\right\rangle$, in other words,

$$\left\langle\mu;U_{1},...,U_{n};\frac{\varepsilon}{2}\right\rangle \subset\left\langle\mu;\varphi;\varepsilon\right\rangle.$$

Let now $\mu=\bigoplus\limits_{i=1}^n\lambda_i\odot\delta_{x_i}\in I_f(X)$, $r_{\delta(X)}^{I(X)}(\mu)=\delta_{x_{i_0}}$, and $U_{i}$ be neighbourhoods of points $x_{i}$ in $X$, respectively, $\left\langle\mu; U_{1}, ..., U_{n}; \varepsilon\right\rangle$ be a neighbourhood of an idempotent probability measure $\mu$ in $I_f(X)$. Denote $V=U_{i_0}$. It is easy to check that at $0<\varepsilon< ln3$ relations $V\cap U_i=\varnothing$, $i=1, ..., n$, $i\neq i_0$ take place. Therefore $r_{\delta(X)}^{I(X)}(\nu)\in \delta(V)$ for any $\nu \in \left\langle\mu; U_{1}, ..., U_{n}; \varepsilon\right\rangle$, where $\delta(V)=\{\delta_x: x\in V\}$ is open in $\delta(X)$ set. From here $r_{\delta(X)}^{I(X)}(\left\langle\mu; U_{1}, ..., U_{n}; \varepsilon\right\rangle) \subset \delta(V)$. Now for each point $y\in V$ we construct an idempotent probability measure
$$\mu_y=0\odot \delta_y\oplus\bigoplus\limits_{\substack{i=1 \\ i\neq i_{0}}}^n\lambda_i\odot\delta_{x_i}.$$
Then, as it is easy to see, $\mu_y\in \left\langle\mu; U_{1}, ..., U_{n}; \varepsilon\right\rangle$ and $r_{\delta(X)}^{I(X)}(\mu_y)=\delta_{y}$. Thus, for each point $y\in V$ we have $\left(r_{\delta(X)}^{I(X)}\right)^{-1}(\delta_y) \in \left\langle\mu; U_{1}, ..., U_{n}; \varepsilon\right\rangle$, i. e. $r_{\delta(X)}^{I(X)}(\left\langle\mu; U_{1}, ..., U_{n}; \varepsilon\right\rangle) = \delta(V)$. Openness of the map $r_{\delta(X)}^{I(X)}$ is established. Theorem 1 is proved.

{\bf Proposition 1.} {\it For any compact $X$  the space $\delta(X)$ is a strong deformation retract of the compact $I_{f}(X)$.}

{\bf Proof.} Consider a map $h:I_{f}(X)\times[0,1]\rightarrow I_{f}(X)$, determined by a formula
$$h(\mu,t)=h_{t}(\mu)=\left(ln(1-t)-lnt\oplus ln(1-t)\right)\odot\mu\oplus\left(lnt-lnt\oplus ln(1-t)\right)\odot r_{f}^{x}(\mu),$$
$$(\mu,t)\in I_{f}(X)\times[0,1].$$
It is easy to check that the map $h$ is defined correctly. Moreover, $h_{0}=id_{I_{f}(X)}$  and  $h_{1}=r_{\delta(X)}^{I(X)}$, i. e.  $h$ is the homotopy connecting maps $id_{I_{f}(X)}$  and  $r_{\delta(X)}^{I(X)}$. Further, we have
$$h(\delta_{x},t)=\left(ln(1-t)-lnt\oplus ln(1-t)\right)\odot\delta_{x}\oplus \left(lnt-lnt\oplus ln(1-t)\right)\odot r_{f}^{x}(\delta_{x})=\delta_{x},$$
i. e. $h_{t}(\delta_{x})=\delta_{x}$  for all $\delta_{x}\in\delta(X)$  and  $t\in[0,1]$. Thus, $\delta(X)$  is a strong deformation retract of a compact $I_{f}(X)$. Proposition 1 is proved.

{\bf Proposition 2.} {\it For any finite compact $X$ the set $I_{f}(X)$  is an neighbourhood retract of the compact $I(X)$.}

{\bf Proof.} For each measure of Dirac $\delta_{x}$, $x\in X$, we build an open set $$O\delta_x=\left\{\mu=\bigoplus\limits_{i=1}^n\lambda_i\odot \delta_{x_i}\in I(X): \lambda_i=0 \mbox{ at } x_i=x \mbox{ and } \lambda_i<-ln 2 \mbox{ at } x_i\neq x \right\}.$$
It is clear, that $O\delta_x\cap O\delta_y=\varnothing$ at $x\neq y$.

Consider open set $\bigcup\limits_{x\in X}O\delta_x$ in $I(X)$. We have $I_{f}(X)\subset\bigcup\limits_{x\in X}O\delta_x$. Now we will show that $I_{f}(X)$ is retract of the open set $\bigcup\limits_{x\in X}O\delta_x$.
Take $\nu\in\bigcup\limits_{x\in X}O\delta_x$, where  $\nu=\bigoplus\limits_{i=1}^{n}\lambda_{i}\odot\delta_{x_{i}}$,  $\bigoplus\limits_{i=1}^{n}\lambda_{i}=\mathbf{1}$,  $\lambda_{i}\geq\mathbf{0}$,  $i=1,...,n$. Then $\lambda_{i_{0}}=0$  for some (only unique) $i_{0}$  and, therefore  $r_{\delta(X)}^{I(X)}(\nu)=\delta_{x_{i_{0}}}$. Define a map $r:\bigcup\limits_{x\in X}O\delta_x\rightarrow I_{f}(X)$  by the rule\\
$
r(\nu) =
$
\[=\begin{cases}
0\odot\delta_{x_{i_{0}}}\oplus\bigoplus\limits_{\substack{i=1,\\ i\neq i_{0}}}^{n} (-ln(n+1))\odot\delta_{x_{i}}, & \text{if  $\lambda_i\geq-ln(n+1)$ for any $i\in \{1, ..., n\}\setminus\{i_0\}$}, \\
\bigoplus\limits_{i=1}^{n}\lambda_{i}\odot\delta_{x_{i}}, & \text{if $\lambda_{i}\leq -ln(n+1)$ for all $i\in \{1, ..., n\}\setminus\{i_0\}$.}
\end{cases}
\]

The map $r$ is defined correctly. It is continuous. Besides $r(\mu)=\mu$  for any measure  $\mu\in I_{f}(X)$. It means the map $r$  is a retraction. Proposition 2 is proved.

Remind [13] that a set $A\subset X$  is called to be collapsible by space $X$ to a set $B\subset X$ if embedding map $i_{A}:A\rightarrow X$ is homotopic to some map $f:A\rightarrow X$ such that $f(A)\subset B$. If $B$  consists of one point, they say that $A$ is collapsible by $X$.

It is clear, if there is a homotype $h:A\times I\rightarrow A$, such that $h(y,0)=i_{A}$, and $h(y,1)=\{\text {a point}\}$ then $A$ collapsible by $X$.

A space $X$  is called to be locally collapsible to a point $x_{0}\in X$ if any neighbourhood $U$ of the point $x_{0}$ contains a neighbourhood $U_{0}$ such, that collapsible to by $U$ to a point. A space $X$ is called to be locally collapsible if it is locally collapsible to each point.

{\bf Theorem 2.} {\it The functor $I_{f}$ preserves a collapsibility of compacts i. e. if $X$ is collapsible compact, then $I_{f}(X)$ is also collapsible compact.}

{\bf Proof.} We will show more: the functor $I_{f}$ preserves a homotopy of maps. Let  $h_{0}$, $h_{1}:X\rightarrow Y$ be homotopic maps,  $h:X\times[0,1]\rightarrow Y$ is the homotopy connecting the maps $h_{0}, \ \ h_{1}$, i. e. $h(x,0)=h_{0}(x)$,  $h(x,1)=h_{1}(x)$. The identity map  $i_{t_{0}}:X\times\{t_{0}\}\rightarrow X\times I$, determined by equality  $i_{t_{0}}(x,t_{0})=(x,t_{0})$,  $x\in X$, defines an identity map $I_{f}(i_{t_{0}}):I_{f}(X\times\{t_{0}\})\rightarrow I_{f}(X\times I)$. But, for every $t_{0}\in[0,1]$ the space  $I_{f}(X\times\{t_{0}\})$ is naturally homeomorphic to $I_{f}(X)\times\{t_{0}\}$. This homeomorphism can be carried out, as it is easy to see, by means of correspondence  $\mu_{t_{0}}\leftrightarrow(\mu,t_{0})$, where $\mu_{t_{0}}=\bigoplus\limits_{i=1}^{n}\lambda_{i}\odot\delta_{(x_{i},t_{0})}\in I_{f}(X\times\{t_{0}\})$  and  $\mu=\bigoplus\limits_{i=1}^{n}\lambda_{i}\odot\delta_{x_{i}}\in I_{f}(X)$.

Define now a map $I_{f}(h):I_{f}(X)\times[0,1]\rightarrow I_{f}(Y)$ by equality  $I_{f}(h)\left(\bigoplus\limits_{i=1}^{n}\lambda_{i}\odot\delta_{x_{i}},t\right)=\bigoplus\limits_{i=1}^{n}\lambda_{i}\odot\delta_{h(x_{i},t)}$. We have
$$I_{f}(h)\left(\bigoplus_{i=1}^{n}\lambda_{i}\odot\delta_{x_{i}},0\right)=\bigoplus_{i=1}^{n}\lambda_{i}\odot\delta_{h(x_{i},0)}=\bigoplus_{i=1}^{n}\lambda_{i}\odot\delta_{h_{0}(x_{i})}=I_{f}(h_{0})\left(\bigoplus_{i=1}^{n}\lambda_{i}\odot\delta_{x_{i}}\right),$$ $$I_{f}(h)\left(\bigoplus_{i=1}^{n}\lambda_{i}\odot\delta_{x_{i}},1\right)=\bigoplus_{i=1}^{n}\lambda_{i}\odot\delta_{h(x_{i},1)}=\bigoplus_{i=1}^{n}\lambda_{i}\odot\delta_{h_{1}(x_{i})}=I_{f}(h_{1})\left(\bigoplus_{i=1}^{n}\lambda_{i}\odot\delta_{x_{i}}\right),$$
i. e. $I_{f}(h)(\mu,0)=I_{f}(h_{0})(\mu)$  and $I_{f}(h)(\mu,1)=I_{f}(h_{1})(\mu)$   for any $\mu\in I_{f}(X)$. In other words, $I_{f}(h)$ is the homotopy connecting the maps $I_{f}(h_{0})$ and $I_{f}(h_{1})$. Thus, functor $I_{f}$ preserves a homotopyness of maps. Theorem 2 is proved.

{\bf Lemma 1.} {\it For a compact $X$ the set $I_{f}(X)$ is the neighbourhood retract of the space  $I_{\omega}(X)$.}

{\bf Proof.} Consider the following set

$$O_{\delta_x}=\left\{\mu=\bigoplus\limits_{i=1}^n \lambda_i\odot\delta_{x_i}\in I_\omega(X):\lambda_i=0 \mbox{ at } x_i=x \mbox{ and } \lambda_i<-ln2 \mbox{ at }x_i\neq x\right\}.$$

We will show that $O_{\delta_x}$ is open in $I_\omega(X)$. Let $\nu=\bigoplus\limits_{i=1}^k \lambda_i\odot\delta_{x_i}\in O_{\delta_x}$. Then $\lambda_i=0$  at $x_i=x$ and $\lambda_i<-ln2$ at $x_i\neq x$. For points $x_1,\ ...,\ x_k$  we choose disjoint open neighbourhoods $U_i$, $i=1,\ ...,\ k$. Put $\varepsilon=\min\left\{-\lambda_i-ln2:\ i\in\{1,\ ...,\ k\}\setminus \{i_0\}\right\} >0$. Let $\mu= \bigoplus\limits_{j=1}^s \gamma_j\odot\delta_{y_j}\in \left\langle\nu;U_1, ... ,U_n;\varepsilon\right\rangle$. Then $supp\mu\cap U_i\neq \varnothing$, and  $supp\mu\subset\bigcup\limits_{i=1}^k U_i$ and $\bigoplus\limits_{y_j\in U_i} |\lambda_i-\gamma_j|<\varepsilon$. On construction for every $i$ and for all $j$, $y_j\in U_i$, we have
$$|\lambda_i-\gamma_j|\leq \bigoplus\limits_{y_j\in U_i} |\lambda_i-\gamma_j|<\varepsilon= \min\left\{-\lambda_i-ln2:\ i=1,\ ...,\ k, i\neq i_0\right\}\leq -\lambda_i-ln2,$$
i. e. $|\lambda_i-\gamma_j|<-\lambda_i-ln2.$  From where we have
$\gamma_j<-ln2$. Hence, $\left\langle\nu;U_1, ... ,U_n;\varepsilon\right\rangle\subset O_{\delta_x}$. Thus, the set $O_{\delta_x}$ is open.

It is obvious that $I_f(X)\subset \bigcup\limits_{x\in X}O_{\delta_x}$. If $\nu\in \bigcup\limits_{x\in X}O_{\delta_x}$, that exists only $x\in X$, that $\nu\in O_{\delta_x}$.Let's construct map $r:\bigcup\limits_{x\in X}O_{\delta_x}\rightarrow I_{f}(X)$  as follows\\

$
r(\nu) =
$
\[=\begin{cases}
0\odot\delta_{x_{i_{0}}}\oplus\bigoplus\limits_{\substack{i=1,\\ i\neq i_{0}}}^{n} (-ln(n+1))\odot\delta_{x_{i}}, & \text{if  $\lambda_i\geq-ln(n+1)$ for some $i\in \{1, ..., n\}\setminus\{i_0\}$}, \\
\bigoplus\limits_{i=1}^{n}\lambda_{i}\odot\delta_{x_{i}}, & \text{if $\lambda_{i}\leq -ln(n+1)$ for all $i\in \{1, ..., n\}\setminus\{i_0\}$,}
\end{cases}
\]
$\nu=\bigoplus\limits_{i=1}^n\lambda_i\odot\delta_{x_i}\in O_{\delta_x}$, $\lambda_{i_0}=0$ for unique $i_0$ and $\lambda_i<-ln2$ at all $i\in\left\{1, ... , n\right\}\backslash \{i_0\}$.

It is easy to see, the map $r$ is a retraction. Lemma 1 is proved.

{\bf Theorem 3.} {\it Let $X$ be an $A(N)R$-compact. Then $I_{f}(X)$ is also $ANR$-compact.}

{\bf Proof.} Let $X$ be neighbourhood retract some compact  $Y$, $U$ be open set in $Y$  such that $U\supset X$  and there is a retraction  $r:U\rightarrow X$.

Consider open set $\bigcup\limits_{x\in U} O\delta_x$ in $I_{f}(Y)$. As we above note $I_{f}(X)\subset \bigcup\limits_{x\in U} O\delta_x$, and if $\nu\in \bigcup\limits_{x\in X}O_{\delta_x}$ then there exists unique $x\in X$ such that $\nu\in O_{\delta_x}$. For  $\nu\in\bigcup\limits_{x\in U} O\delta_x\subset I_{f}(Y)$,  $\nu=\bigoplus\limits_{i=1}^{n}\lambda_{i}\odot\delta_{y_{i}}$,  $\bigoplus\limits_{i=1}^{n}\lambda_{i}=\mathbf{1}$,  $\lambda_{i}\geq \mathbf{0}$,  $i=1,2,...,n$,  $\lambda_{i_{0}}= 0$ put
 $$r_{U}^{Y}(\nu)=\left(\lambda_{i_{0}}\oplus\bigoplus\limits_{y_{i}\in Y\backslash U}\lambda_{i}\right)\odot\delta_{y_{i_{0}}}\oplus\bigoplus\limits_{y_{i}\in U}\lambda_{i}\odot\delta_{y_{i}}.$$
It is obvious that $r_{U}^{Y}(\nu)\in \bigcup\limits_{x\in U} O\delta_x$. Besides, $r_{U}^{Y}(\nu)=\nu$  for any measure  $\nu\in I_{f}(Y)$ such that  $supp\nu\subset U$. As the operation of capture of a maximum is continuous, the constructed map $r_{U}^{Y}:\bigcup\limits_{x\in U} O\delta_x\rightarrow \bigcup\limits_{x\in U} O\delta_x$ is continuous. Further, we will put
$$R\left(r_{U}^{Y}(\nu)\right)=\left(\lambda_{i_{0}}\oplus\bigoplus\limits_{y_{i}\in Y\backslash U}\lambda_{i}\right)\odot\delta_{r\left(y_{i_{0}}\right)}\oplus\bigoplus\limits_{y_{i}\in U}\lambda_{i}\odot\delta_{r\left(y_{i}\right)}.$$
On construction we have $R\left(r_{U}^{Y}(\nu)\right)\in I_{f}(X)$. The map $R:\bigcup\limits_{x\in U} O\delta_x\rightarrow I_{f}(X)$  is defined correctly. As retraction $r:U\rightarrow X$  is continuous, the map $R$ is also continuous. It is easy to check that $R\left(r_{U}^{Y}(\nu)\right)=\nu$  for any measure  $\nu\in I_{f}(X)$. Thus, $R\circ r_{U}^{Y}:\bigcup\limits_{x\in U} O\delta_x\rightarrow I_{f}(X)$  is required retraction. So, the set $I_{f}(X)$ is neighbourhood retract of the compact  $I_{f}(Y)$.

Now application of a Lemma 1 and theorem (3.1) from work [13] finishes the proof of Theorem 3.

In funktorial mean Theorem 3 looks as so:

{\bf Corollary 1.} {\it Functor  $I_{f}$ preserves $ANR$-compacts.}

Theorem 3 and Proposition 1 imply the following important result.

{\bf Corollary 2.} {\it Let $X$  be a compact. $I_{f}(X)\in ANR$  if and only if  $X\in ANR$.}

Further, Theorem 3 implies the following statements.

{\bf Corollary 3.} {\it Functor $I_{f}$  preserves property of a compacts to be $Q$-manifold or Hilbert cube.}

{\bf Corollary 4.} {\it Functor $I_{f}$ preserves property of fibres of maps to be $ANR$-compact, compact  $Q$-manifold and Hilbert cube (finite sum of Hilbert cube).}

\begin{center}
\textsl{References}
\end{center}

{1.} G.L.Litvinov, V.P.Maslov(eds.), Idempotent mathematics and mathematical physics (Vienna, 2003), Contemp.Math.,377,Amer.Math.Soc., Providence, RI, 2005.

{2.} G.L.Litvinov, ``The Maslov dequantization, idempotent and tropical mathematics: a very brief introduction''
Idempotent mathematics and mathematical physics (Vienna, 2003), Contemp. Math., 377, Amer.Math.Soc., Providence, RI, 2005, 1-17; arXiv:abs/math/0501038.

{3.} P.Bernhard, ``Max-plus algebra and mathematical fear in dynamic optimization'', set-Valued Anal., 8:1-2 (2000), 71-84.

{4.} J.P.Aubin, O.Dordan, ``Fuzzy systems, viability theory and toll sets'', Fuzzy systems, Handb.Fuzzy Sets Ser., 2, Kluwer Acad.Publ., Boston, MA, 1998, 461-488.

{5.} J.-P. Aubin, Dynamic economic theory. A viability approach, Stdu. Econom. Theory, 5, Springer-Verlag, Berlin, 1997.

{6.} V.V.Fedorchuk, Probability measures in topolog'', UMN, 46:1 (1981), 41-80; trans.engl.: Russian Math. Surveys, 46:1 (1991), 45-93.

{7.} M. Zarichnyi. Idempotent probability measures,
I.arXiv:math.GN/0608754 v 130 Aug 2006.

{8.} Al-Kassas Yusef. Metrizability and paracompactness of the space of probability measures. PhD thesis. Moscow state university. 1991 (Rusian).

{9.} E. V. Shepin. Functors and incalculable degrees of compacts. UMN. 1981,
т36, extract 3, p.3-62.

{10.} O. H. Keller. Die Homeomorphieder Kompakten convexen Mengenin Hilbertschen Raum. Math.Ann. 1931, Bd 105, 8, 748-758

{11.} T. A. Chapman. Lectures on Hilbert Cube Manifolds. American Mathematical Society, Providence, Rhode Island, National Science Foundation. 1975.

{12.} K. Borsuk. Theory of shape.  Publisher: Taylor and Francis. June 1 1975.

{13.} K. Borsuk. Theory of retracts. Warszawa. 1967. PWN -- Polish scientific publishers.

\end{document}